\documentclass[12pt]{report}
\usepackage{amssymb,latexsym,epsfig,mathrsfs,graphpap,graphics,amsmath,amssymb}
\usepackage{amstext}
\usepackage{amsthm}
\usepackage[utf8]{inputenc}
\usepackage[english]{babel}
\usepackage{helvet}
\usepackage{courier}

\theoremstyle{Definition}

\theoremstyle{remark}

\numberwithin{equation}{section}
\begin{document}

\begin{quote}

 {\bf\Large {Construction of $J^{\text{th}}$-stage Nonuniform   Wavelets on Local Fields }}

\parindent=0mm \vspace{.3in}

  {\bf{Owais Ahmad$^{\star}$ and Firdous A. Shah$^{\star\star}$}}

  \parindent=0mm \vspace{.1in}
{\it\small{$^{\star}${Department of Mathematics, National Institute of Technology, Srinagar-190 006, Jammu and Kashmir, India. E-mail: $\text{siawoahmad@gmail.com}$}}

$^{\star\star}${Department of  Mathematics,  University of Kashmir, South Campus, Anantnag-192 101, Jammu and Kashmir, India. E-mail: $\text{fashah79@kashmiruniversity.ac.in}$ }}

\parindent=0mm \vspace{.1in}
 {\small {\bf Abstract:}  Shah and Abdullah [Complex Analysis Operator Theory, 9 (2015), 1589-1608] have introduced  a generalized notion of nonuniform multiresolution analysis (NUMRA) on local field $K$ of positive characteristic  in which the translation set $\Lambda$ acting on the scaling function to generate the core space $V_{0}$ is no longer a group, but is the union of ${\mathcal Z}$ and a translate of ${\mathcal Z}$, given by
 $\Lambda=\left\{0,u(r)/N \right\}+{\mathcal Z}$, where $N \ge 1$ is an integer and $r$ is an odd integer such that $r$ and $N$ are relatively prime, and ${\mathcal Z}=\{u(n): n\in\mathbb N_{0}\}$  is a complete list of distinct cosets  of the unit disc $\mathfrak D$ in $K^+.$  In this paper,
 we focus on the extension of nonuniform continuous wavelets to the construction of $J^{\text{th}}$-stage nonuniform discrete wavelets on local fields. We  establish some general characterizations for the $J^{\text{th}}$-stage nonuniform discrete wavelet systems to be  orthornormal bases in $L^2(\Lambda)$. Moreover, we establish a relation between the continuous wavelets of $L^2(K)$ and their discrete counterparts of $l^2(\Lambda)$.


\parindent=0mm \vspace{.0in}
 }

\end{quote}

\parindent=0mm \vspace{.2in}
{\bf{1. Introduction}}

\parindent=0mm \vspace{.1in}
In recent years, there has been a considerable interest in the study of harmonic analysis and wavelet analysis over the local fields. Local fields are essentially of two types: zero and positive characteristic. Examples of local fields of characteristic zero include the $p$-adic field $\mathbb Q_{p}$ where as local fields of positive characteristic are the Cantor dyadic group and the Vilenkin $p$-groups.  
Despite the fact that the structures and metrics of  $p$-adic fields and local fields of positive characteristic are comparable, their wavelet and MRA theory are quite different. The notion of  MRA on  local fields of positive characteristic was introduced by Jiang et al.\cite{JLJ}. In fact, they brought up a technique for constructing orthogonal wavelets on  local fields  and established a necessary and sufficient condition for the solution of  refinement equation to generate an MRA for $L^2(K)$. Subsequently, an explicit construction of tight wavelet frames on   local fields was given by Shah and Debnath \cite{Shahdb} by adapting the  extension principles   on the Euclidean spaces to the local fields. On the other hand, Shah and Abdullah \cite{ShahA1} have set up an entire portrayal of tight wavelet frames on local fields by virtue of some fundamental equations in the frequency domain and demonstrate how to build the Parseval wavelet frames for $L^2(K)$. These studies were proceeded by Shah and his associates in \cite{Shah1,ShahA2,ShahO1,ShahY2}, where they have given some algorithms for constructing periodic wavelet frames, wave packet frames,  and semi-orthogonal wavelet frames on non-Archimedean local fields of positive characteristic.

\parindent=8mm \vspace{.2in}
In our previous work \cite{ShahA3}, we  have generalized the concept of Mallat's classic MRA  on Euclidean spaces $\mathbb R^n$ to  nonuniform MRA on local fields of positive characteristic, in which the translation set acting on the scaling function associated with the MRA to generate the core space $V_{0}$ is no longer a group, but is the union of ${\mathcal Z}$ and a translate of ${\mathcal Z}$, where ${\mathcal Z}=\{u(n): n\in\mathbb N_{0}\}$  is a complete list of (distinct) coset representation of the unit disc $\mathfrak D$ in the locally compact Abelian group $K^+.$  More precisely, this set is of the form $\Lambda=\left\{0,u(r)/N \right\}+{\mathcal Z}$, where $N \ge 1$ is an integer and $r$ is an odd integer such that $r$ and $N$ are relatively prime. We  call this a {\it  nonuniform multiresolution analysis} (NUMRA) on local fields of positive characteristic. As a consequence of this generalization, we obtain a necessary and sufficient condition for the existence of associated wavelets and extension of Cohen's theorem. Recently, we have  constructed the associated  nonuniform wavelet packets  on local fields in \cite{ShahY3}. Indeed, we obtain a lemma on the so-called splitting trick and several theorems concerning the Fourier transform of the nonuniform wavelet packets  to show that their translates form an orthonormal basis for $L^2(K)$. More results in this direction can also be found in \cite{ShukV} and the references therein.

\parindent=8mm \vspace{.2in}
Owing to the fact that the data in both physics and engineering is often discrete in nature, makes us to focus our investigation over the discrete sequence spaces on local fields of positive characteristic. The concept of an adaptive MRA structure was introduced by Han et al.\cite{han2} for more general affine-like systems which exhibits all the favorable properties of MRA structures for wavelets whereas Han \cite{han1} has independently developed a comprehensive theory of discrete framelets and wavelets using an algorithmic approach by directly studying a discrete framelet transform. The main contribution of this paper is that we extend our previous work \cite{ShahA3} and construct a class of $J^{\text{th}}$-stage nonuniform discrete wavelet systems on local fields of positive characteristic. Different from our previous approach in the orthonormal case, our analysis of nonuniform discrete scheme is inspired by Shukla and Mittal's approach in \cite{ShukMV} for construction of  wavelets on the spectrum. We  provide some characterizations of the $J^{\text{th}}$-stage discrete wavelet systems to be orthonormal bases for the Hilbert space  $\l^2(\Lambda)$. Moreover, we establish a connection between a system of nonuniform wavelets of $L^2(K)$ and a first-stage nonuniform discrete wavelet system of $l^2(\Lambda)$.

\pagestyle{myheadings}

\parindent=8mm \vspace{.2in}
The article is organized as follows. In Sect. 2, we give a necessary background about local  fields including the definitions of Fourier transform, uniform MRA and non-uniform MRA on fields fields. In Sect. 3, we introduce the construction of a first-stage nonuniform discrete wavelet system and provide a characterization  for such a system to be an orthonormal basis  for the Hilbert space $l^2(\Lambda)$. Sect. 4, is devoted to the construction of $J^{\text{th}}$-stage nonuniform discrete wavelets for $l^2(\Lambda)$ by decomposing of its closed subspaces. Finally,  we establish a relation between the continuous wavelets of $L^2(K)$ and their discrete counterparts of $l^2(\Lambda)$ in Sect. 5.

\parindent=0mm \vspace{.2in}
{\bf{2. Fourier and Wavelet Analysis on  local  Fields }}

\parindent=0mm \vspace{.2in}
In this section, we present some important preliminaries and notation that will be useful in the sequel to obtain certain characterizations of $J^{\text{th}}$-stage nonuniform discrete orthonormal wavelet bases for $l^2(\Lambda)$. More precisely, we review some concepts about Fourier and wavelet analysis on local fields of positive characteristic.

\parindent=0mm \vspace{.1in}
{\it {2.1. local Fields}}

\parindent=0mm \vspace{.1in}
A  local field $K$ is a  locally compact, non-discrete and totally disconnected field. If it is of characteristic zero, then  it is a field of $p$-adic numbers $\mathbb Q_p$ or its finite extension. If $K$ is of positive characteristic, then $K$ is a field of formal Laurent series over a finite field $GF(p^c)$. If $c =1$, it is a $p$-series field, while for $c\ne 1$, it is an algebraic extension of degree  $c$ of a $p$-series field. Let $K$ be a fixed local  field  with the ring of integers ${\mathfrak D}= \left\{x \in K: |x| \le 1\right\}$. Since $K^{+}$ is a locally compact Abelian group, we choose a Haar measure $dx$ for $K^{+}$. The  field $K$ is locally compact, non-trivial, totally disconnected and complete topological field endowed with non--Archimedean norm  $|\cdot|:K\to \mathbb R^+$ satisfying

\parindent=0mm \vspace{.1in}
(a) $|x|=0$ if and only if $x = 0;$

\parindent=0mm \vspace{.1in}
(b) $|x\,y|=|x||y|$ for all $x, y\in K$;

\parindent=0mm \vspace{.1in}
(c) $|x+y|\le \max \left\{ |x|, |y|\right\}$ for all $x, y\in K$.

\parindent=0mm \vspace{.1in}
Property (c) is called the ultrametric inequality. Let ${\mathfrak B}= \left\{x \in K: |x| < 1\right\}$ be the prime ideal of the ring of integers ${\mathfrak D}$ in $K$. Then, the residue space ${\mathfrak D}/{\mathfrak B}$ is isomorphic to a finite field $GF(q)$, where $q = p^{c}$ for some prime $p$ and $c\in\mathbb N$. Since  $K$ is totally disconnected and $\mathfrak B$ is both prime and principal ideal, so there exist a prime element $\mathfrak p$ of $K$ such that ${\mathfrak B}= \langle \mathfrak p \rangle=\mathfrak p {\mathfrak D}$. Let ${\mathfrak D}^*= {\mathfrak D}\setminus {\mathfrak B }=\left\{x\in K: |x|=1   \right\}$. Clearly,  ${\mathfrak D}^*$ is a group of units in $K^*$ and if $x\not=0$, then can write $x=\mathfrak p^n y, y\in {\mathfrak D}^*.$ Moreover, if ${\cal U}= \left\{a_m:m=0,1,\dots,q-1 \right\}$ denotes the fixed full set of coset representatives of ${\mathfrak B}$ in ${\mathfrak D}$, then every element $x\in K$ can be expressed uniquely  as $x=\sum_{\ell=k}^{\infty} c_\ell \,\mathfrak p^\ell $ with $c_\ell \in {\cal U}.$ Recall that ${\mathfrak B}$ is compact and open, so each  fractional ideal ${\mathfrak B}^k= \mathfrak p^k {\mathfrak D}=\left\{x \in K: |x| < q^{-k}\right\}$  is also compact and open and is a subgroup of $K^+$. We use the notation in Taibleson's book \cite{Tab}. In the rest of this paper, we use the symbols $\mathbb N, \mathbb N_0$ and $\mathbb Z$ to denote the sets of natural, non-negative integers and integers, respectively.

\parindent=8mm \vspace{.1in}
 Let $\chi$ be a fixed character on $K^+$ that is trivial on ${\mathfrak D}$ but  non-trivial on  ${\mathfrak B}^{-1}$. Therefore, $\chi$ is constant on cosets of ${\mathfrak D}$ so if $y \in {\mathfrak B}^k$, then $\chi_y(x)=\chi(y,x), x\in K.$ Suppose that $\chi_u$ is any character on $K^+$, then the restriction $\chi_u|{\mathfrak D}$ is a character on ${\mathfrak D}$. Moreover, as characters on ${\mathfrak D}, \chi_u=\chi_v$ if and only if $u-v\in {\mathfrak D}$. Hence, if  $\left\{u(n): n\in\mathbb N_0\right\}$ is a complete list of distinct coset representative of ${\mathfrak D}$ in $K^+$, then, as it was proved in \cite{Tab}, the set  $\left\{\chi_{u(n)}: n\in\mathbb N_0\right\}$   of distinct characters on ${\mathfrak D}$ is a complete orthonormal system on ${\mathfrak D}$.

\parindent=8mm \vspace{.1in}
We now impose a natural order on the sequence $\{u(n)\}_{n=0}^\infty$. We have ${\mathfrak D}/ \mathfrak B \cong GF(q) $ where $GF(q)$ is a $c$-dimensional vector space over the field $GF(p)$. We choose a set $\left\{1=\zeta_0,\zeta_1,\zeta_2,\dots,\zeta_{c-1}\right\}\subset {\mathfrak D^*}$ such that span$\left\{\zeta_j\right\}_{j=0}^{c-1}\cong GF(q)$. For $n \in \mathbb N_0$ satisfying
$$0\leq n<q,~~n=a_0+a_1p+\dots+a_{c-1}p^{c-1},~~0\leq a_k<p,~~\text{and}~k=0,1,\dots,c-1,$$

\parindent=0mm \vspace{.1in}
we define
$$u(n)=\left(a_0+a_1\zeta_1+\dots+a_{c-1}\zeta_{c-1}\right){\mathfrak p}^{-1}.\eqno(2.1)$$

\parindent=0mm \vspace{.1in}
Also, for $n=b_0+b_1q+b_2q^2+\dots+b_sq^s, ~n\in \mathbb N_{0},~0\leq b_k<q,k=0,1,2,\dots,s$, we set

$$u(n)=u(b_0)+u(b_1){\mathfrak p}^{-1}+\dots+u(b_s){\mathfrak p}^{-s}.\eqno(2.2)$$

\parindent=0mm \vspace{.1in}
This defines $u(n)$ for all $n\in \mathbb N_{0}$. In general, it is not true that $u(m + n)=u(m)+u(n)$. But, if $r,k\in\mathbb N_{0}\; \text{and}\;0\le s<q^k$, then $u(rq^k+s)=u(r){\mathfrak p}^{-k}+u(s).$ Further, it is also easy to verify that $u(n)=0$ if and only if $n=0$ and $\{u(\ell)+u(k):k \in \mathbb N_0\}=\{u(k):k \in \mathbb N_0\}$ for a fixed $\ell \in \mathbb N_0.$ Hereafter we use the notation $\chi_n=\chi_{u(n)}, \, n\ge 0$.

\parindent=8mm \vspace{.2in}
Let the local field $K$ be of characteristic $p>0$ and $\zeta_0,\zeta_1,\zeta_2,\dots,\zeta_{c-1}$ be as above. We define a character $\chi$ on $K$ as follows:
$$\chi(\zeta_\mu {\mathfrak p}^{-j})= \left\{
\begin{array}{lcl}
\exp(2\pi i/p),&&\mu=0\;\text{and}\;j=1,\\
1,&&\mu=1,\dots,c-1\;\text{or}\;j \neq 1.
\end{array}
\right. \eqno(2.3)$$

\parindent=0mm \vspace{.1in}
{\it {2.2.  Fourier Transforms on local Fields}}

\parindent=0mm \vspace{.2in}
 The Fourier transform of $f \in L^1(K)$ is denoted by $\hat f(\xi)$ and defined  by
\begin{align*}
{\mathscr F}\big\{f(x)\big\}=\hat f(\xi)=\int_K f(x)\overline{ \chi_\xi(x)}\,dx.\tag{2.4}
\end{align*}
It is noted that
$$\hat f(\xi)= \displaystyle \int_K f(x)\,\overline{ \chi_\xi(x)}dx= \displaystyle \int_K f(x)\chi(-\xi x)\,dx.$$

\parindent=0mm \vspace{.1in}
The properties of Fourier transforms on local field $K$ are much similar to those of on the classical field $\mathbb R$. In fact, the Fourier transform on local fields of positive characteristic have the following properties:
\begin{itemize}
  \item The map $f\to \hat f$ is a bounded linear transformation of $L^1(K)$ into $L^\infty(K)$, and $\big\|\hat f\big\|_{\infty}\le \big\|f\big\|_{1}$.
  \item If $f\in L^1(K)$, then $\hat f$ is uniformly continuous.
  \item If $f\in L^1(K)\cap L^2(K)$, then $\big\|\hat f\big\|_{2}=\big\|f\big\|_{2}$.
\end{itemize}

\parindent=0mm \vspace{.1in}
The Fourier transform of a function $f\in L^2(K)$ is defined by
\begin{align*}
\hat f(\xi)= \lim_{k\to \infty} \hat f_{k}(\xi)=\lim_{k\to \infty}\int_{|x|\le q^{k}} f(x)\overline{ \chi_\xi(x)}\,dx,\tag{2.5}
\end{align*}

\parindent=0mm \vspace{.0in}
where  $f_{k}=f\,\Phi_{-k}$ and $\Phi_{k}$ is  the characteristic function of ${\mathfrak B}^{k}$.  Furthermore, if $f\in L^2(\mathfrak D)$, then we define the Fourier coefficients of $f$ as
\begin{align*}
\hat f\big(u(n)\big)=\int_{\mathfrak D} f(x) \overline{ \chi_{u(n)}(x)}\,dx.\tag{2.6}
\end{align*}

\parindent=0mm \vspace{.0in}
The series $\sum_{n\in \mathbb N_{0}} \hat f\big( u(n)\big) \chi_{u(n)}(x)$ is called the Fourier series of $f$. From the standard $L^2$-theory for compact Abelian groups, we conclude that the Fourier series of $f$ converges to $f$ in $L^2(\mathfrak D)$ and Parseval's identity holds:
\begin{align*}
\big\|f\big\|^2_{2}=\int_{\mathfrak D}\big|f(x)\big|^2 dx= \sum_{n\in \mathbb N_{0}} \left| \hat f\big(u(n)\big)\right|^2.\tag{2.7}
\end{align*}

\parindent=0mm \vspace{.0in}
{\it {2.3. Uniform MRA on local Fields}}

\parindent=0mm \vspace{.1in}
In order to able to define the concepts of uniform MRA and wavelets on local fields, we need analogous notions of translation and dilation. Since $\bigcup_{j\in\mathbb Z} \mathfrak p^{-j} {\mathfrak D}=K$, we can regard $\mathfrak p^{-1}$ as the dilation and since $\left\{u(n):n\in\mathbb N_{0}\right\}$ is a complete list of distinct coset representatives of $\mathfrak D$ in $K$, the set ${\mathcal Z}=\left\{u(n):n\in\mathbb N_{0}\right\}$ can be treated as the translation set. Note that $\Lambda$ is a subgroup of $K^{+}$ and unlike the standard wavelet theory on the real line, the translation set is not a group.

\parindent=8mm \vspace{.1in}
The following is a definition of uniform MRA on local  fields of positive characteristic \cite{JLJ}.

\parindent=0mm \vspace{.1in}
{\bf{Definition 2.1.}} Let $K$ be a local field of positive characteristic $p>0$ and $ \mathfrak p$ be a prime element of $K$. An MRA of $L^2(K)$ is a  sequence of closed subspaces $\{V_j:j\in \mathbb Z\}$ of $L^2(K)$ satisfying the following properties:

\parindent=0mm \vspace{.1in}
(a)\quad $V_j \subset V_{j+1}\; \text{for all}\; j \in \mathbb Z;$

\parindent=0mm \vspace{.1in}
(b)\quad $\bigcup_{j\in \mathbb Z}V_j\;\text{is dense in}\;L^2(K);$

\parindent=0mm \vspace{.1in}
(c)\quad $\bigcap_{j\in \mathbb Z}V_j=\{0\};$

\parindent=0mm \vspace{.1in}
(d)\quad $f(x) \in V_j\; \text{if and only if}\;f({\mathfrak p}^{-1}x) \in V_{j+1}\; \text{for all}\; j \in \mathbb Z;$

\parindent=0mm \vspace{.1in}
(e) There exists a function $\phi \in V_0$, such that $\left\{\phi\big(x-u(k)\big): k\in \mathbb N_0\right\}$ forms an orthonormal basis for $V_0$.

\parindent=5mm \vspace{.1in}
According to the standard scheme for construction of MRA-based wavelets, for each $j$, we define a wavelet space $W_{j}$ as the orthogonal complement of $V_{j}$ in $V_{j+1}$, i.e., $V_{j+1}=V_{j}\oplus W_{j}, \, j\in\mathbb Z$, where $W_{j}\perp V_{j}, \, j\in\mathbb Z$. It is not difficult to see that
$$f(x)\in W_{j} \quad \text{if and only if}\quad f(\mathfrak p^{-1}x)\in W_{j+1},\quad j\in\mathbb Z.\eqno(2.7)$$

\parindent=0mm \vspace{.0in}
Moreover, they are mutually orthogonal, and we have the following orthogonal decompositions:
$$L^2(K)= \bigoplus_{j\in\mathbb Z} W_{j}=V_{0}\oplus \left(\bigoplus_{j\ge 0}W_{j}\right).\eqno(2.8)$$

\parindent=0mm \vspace{.0in}
As in the case of $\mathbb R^n$, we expect the existence of $q-1$ number of functions $\psi_{1}, \psi_{2},\dots, \psi_{q-1}$ to form a set of basic wavelets. In view of (2.7) and (2.8), it is clear that if $\left\{\psi_{1}, \psi_{2},\dots, \psi_{q-1}\right\}$ is a set of function such that the system $\left\{\psi_{\ell}\big(x-u(k)\big): 1\le \ell\le q-1, k\in\mathbb N_{0} \right\}$ forms an orthonormal basis for $W_{0}$, then  $\left\{q^{j/2}\psi_{\ell}(\mathfrak p ^{-j}x-u(k)\big): 1\le \ell\le q-1,j\in\mathbb Z, k\in\mathbb N_{0} \right\}$ forms an orthonormal basis for $L^2(K)$.

\parindent=0mm \vspace{.2in}
{\it {2.4. Nonuniform MRA  on local Fields}}

\parindent=0mm \vspace{.1in}
 For an integer $N \ge 1$ and an odd integer $r$ with $1\leq r \leq qN-1$ such that $r$ and $N$ are relatively prime, we define $$\Lambda = \left\{ 0, \dfrac{u(r)}{N}\right\}+{\mathcal Z}.$$

 \parindent=0mm \vspace{.1in}
where ${\mathcal Z}=\left\{ u(n): n\in \mathbb N_{0}\right\}$. It is easy to verify that $\Lambda$ is not a group on local field $K$, but is the union of ${\mathcal Z}$ and a translate of ${\mathcal Z}.$ Following is the definition of nonuniform multiresolution analysis (NUMRA) on local fields of positive characteristic given by Shah and Abdullah \cite{ShahA3}.

\parindent=0mm \vspace{.1in}
{\bf{Definition 2.2.}} For an integer $N \ge 1$ and an odd integer $r$ with $1\leq r \leq qN-1$ such that $r$ and $N$ are relatively prime, an associated  NUMRA on local field $K$ of positive characteristic is a sequence of closed subspaces $\left\{V_j: j\in\mathbb Z\right\}$ of $L^2(K)$ such that the following properties hold:

\parindent=0mm \vspace{.2in}
(a)\quad $V_j \subset V_{j+1}\; \text{for all}\; j \in \mathbb Z;$

\parindent=0mm \vspace{.1in}
(b)\quad $\bigcup_{j\in \mathbb Z}V_j\;\text{is dense in}\;L^2(K);$

\parindent=0mm \vspace{.1in}
(c)\quad $\bigcap_{j\in \mathbb Z}V_j=\{0\};$

\parindent=0mm \vspace{.1in}
(d)\quad $f(x) \in V_j\; \text{if and only if}\;f({\mathfrak p}^{-1}Nx) \in V_{j+1}\; \text{for all}\; j \in \mathbb Z;$

\parindent=0mm\vspace{.1in}

(e)~ There exist a function $\phi$ in $V_0$ such that the collection $\left\{ \phi (x- \lambda ): \lambda \in \Lambda\right\}$ is a complete orthonormal basis for $V_0$.

\parindent=8mm \vspace{.1in}
It is worth noticing that, when $N = 1$, one recovers from the definition above the  definition of an MRA on local fields of positive characteristic $p>0$. When, $N > 1$, the dilation is induced by $\mathfrak p^{-1}N$ and $|\mathfrak p^{-1}|=q$ ensures that $qN\Lambda\subset {\mathcal Z} \subset \Lambda$.

\parindent=8mm \vspace{.2in}
As in the standard scheme, one expects the existence of $qN -1$ number of functions so that their translation by elements of $\Lambda$ and dilations by the integral powers of ${\mathfrak p^{-1}}N$ form an orthonormal basis for $L^2(K)$.

\parindent=0mm \vspace{.1in}
{\bf{Definition 2.3.}} A set of functions $\left\{\psi_{1}, \psi_{1},\dots,\psi_{qN-1}\right\}$ in $L^2(K)$ is said to be a {\it set of  basic wavelets} associated with an NUMRA $\left\{V_{j}: j\in\mathbb Z\right\}$ if the family of functions $\left\{ (qN)^{j/2} \psi_{\ell}\left(\big(\mathfrak p^{-1}N\big)^{j}x-\lambda\right):1\le \ell\le qN-1, \lambda\in \Lambda\right\}$ forms an orthonormal basis for $W_j$.

\parindent=0mm \vspace{.2in}
{\bf{3. First-stage Discrete Wavelets on local Fields}}

\parindent=0mm \vspace{.2in}
The main content of this section is to establish a characterization of the first-stage nonuniform discrete wavelets on local fields of positive characteristic.

\parindent=8mm \vspace{.1in}
We regard $z$ as a function defined on the set  $\Lambda$ and suppose that $z=\left\{z(\lambda)\right\}_{\lambda\in \Lambda}$. We define the spaces
\begin{align*}
l^2(\Lambda) &=\Big\{z :\Lambda \to \mathbb C:  \sum_{\lambda \in \Lambda}\big|z(\lambda)\big|^2 < \infty \Big\},\qquad \text{and}\\
L^2(\Omega) &=\left\{f :\Omega \to \mathbb C: \int_\Omega\big|f(\xi)\big|^2 d\xi< \infty \right\},
\end{align*}
where $\Omega$ is a Lebesgue measurable subset of $K$ with finite positive measure. These spaces are Hilbert spaces with the inner products defined by
\begin{align*}
\langle z,w\rangle  &= \sum_{\lambda \in \Lambda}z(\lambda)\overline{w(\lambda)}\quad \text{for}\; z, w \in l^2(\Lambda),\qquad \text{and}\\
\langle f, g\rangle &= \int_\Omega f(\xi)\overline{g(\xi)}\, d\xi\quad \text{for}\; f, g \in L^2(\Omega),
\end{align*}
respectively.

\parindent=0mm \vspace{.1in}
{\bf{Definition 3.1.}} The Fourier transform on $l^2(\Lambda)$ is a map $ \wedge:l^2(\Lambda) \to L^2(\Omega)$ defined by
\begin{align*}
\hat z(\xi)=\sum_{\lambda \in \Lambda}z(\lambda)\overline{\chi_{\lambda}(\xi)}, \quad z \in l^2(\Lambda)\tag{3.1}
\end{align*}
and its inverse is given by
\begin{align*}
f^{\vee}(\lambda)=\left\langle f,\overline{\chi_{\lambda}(\xi)}\right\rangle=\int_{\Omega}f(\xi)\chi_{\lambda}(\xi)\,d\xi , \quad f \in L^2(\Omega).\tag{3.2}
\end{align*}

\parindent=0mm \vspace{.0in}
For all $z,w \in l^2(\Lambda)$,  the Parseval and Plancherel formulae are  given by
\begin{align*}
\langle z,w\rangle = \sum_{\lambda \in \Lambda}z(\lambda)\overline{w(\lambda)}=\langle\hat z, \hat w\rangle, \quad\text{and}\quad
\big|\big|z\big|\big|^2=\sum_{\lambda \in \Lambda}\big|z(\lambda)\big|^2=\big|\big|\hat z\big|\big|^2.
\end{align*}
 For $N\ge 1, q^{-1}=|\mathfrak p|$ and $\lambda \in \Lambda$, the translation operator $T_{qN\lambda}:l^2(\Lambda) \to l^2(\Lambda)$ is defined by
\begin{align*}
T_{qN\lambda}z(\sigma)=z(\sigma-qN\lambda), \quad \forall~\sigma \in \Lambda.
\end{align*}

\parindent=0mm \vspace{.0in}
Then,  for $z,w \in l^2(\Lambda)$, it can be easily verified that
\begin{align*}
\big(T_{qN\lambda}z\big)^{\wedge}(\xi)=\overline{\chi_{qN\lambda}(\xi)}\hat z(\xi)\quad \text{and}\quad \left\langle\, T_{qN\lambda}z, T_{qN\sigma}w\right\rangle=\left\langle T_{qN(\lambda-\sigma)}z, w\right\rangle.\end{align*}

\parindent=0mm \vspace{.1in}
{\bf{Definition 3.2.}}   Let $N \in\mathbb N$ and let $k$ be an odd integer with $1 \le k \le qN-1$ such that $k$ and $N$ are relatively prime. For $w_k\in l^2(\Lambda)$, we call $\mathscr F(W)$ a first stage nonuniform discrete wavelet system associated with $W=\left\{ w_k: w_k\in l^2(\Lambda)\right\}$ if
\begin{align*}
{\mathscr F(W)}=\Big\{T_{qN\lambda}w_k: w_k\in l^2(\Lambda); \lambda \in \Lambda ; 0 \le k \le qN-1\Big\}. \tag{3.3}
\end{align*}
 is a complete orthonormal set in $l^2(\Lambda)$. We shall call $w_0$ as the {\it nonuniform father wavelet} and $\big\{w_k: 1\le k \le qN-1\big\}$ as the {\it nonuniform mother wavelets}.

\parindent=0mm \vspace{.1in}
{\bf{Theorem 3.3.}} {\it For $z,w \in l^2(\Lambda)$, the systems $\left\{ T_{qN\lambda}z\right\}_{\lambda \in \Lambda}$ and $\left\{ T_{qN\lambda}w\right\}_{\lambda \in \Lambda}$ generates orthogonal subspaces in $l^2(\Lambda)$ if and only if the following conditions hold:}
\begin{align*}
&\text{(a)}~\sum_{s=0}^{qN-1}\left\{\hat z\left(\xi+\dfrac{ u(s)}{\mathfrak p^{-1}N}\right)\overline{\hat w\left(\xi+\dfrac{ u(s)}{\mathfrak p^{-1}N}\right)}\hat z\left(\xi+\dfrac{ u(s)}{\mathfrak p^{-1}N}+u(N)\right)\overline{\hat w\left(\xi+\dfrac{ u(s)}{\mathfrak p^{-1}N}+u(N)\right)}\right\}=0,\\
&\text{(b)}~\sum_{s=0}^{qN-1}\overline{\chi\left(\dfrac{r}{N}\, \mathfrak p u(s) \right)}\left\{\hat z\left(\xi+\dfrac{ u(s)}{\mathfrak p^{-1}N}\right)\overline{\hat w\left(\xi+\dfrac{ u(s)}{\mathfrak p^{-1}N}\right)}\right.\\
&\qquad\qquad\qquad\qquad \qquad\qquad  \left.\hat z\left(\xi+\dfrac{ u(s)}{\mathfrak p^{-1}N}+u(N)\right)\overline{\hat w\left(\xi+\dfrac{ u(s)}{\mathfrak p^{-1}N}+u(N)\right)}\right\}=0, \\
&\text{(c)}~\sum_{s=0}^{qN-1}\chi\left(\dfrac{r}{N}\, \mathfrak p u(s) \right)\left\{\hat z\left(\xi+\dfrac{ u(s)}{\mathfrak p^{-1}N}\right)\overline{\hat w\left(\xi+\dfrac{ u(s)}{\mathfrak p^{-1}N}\right)}\right.\\
&\qquad\qquad\qquad\qquad \qquad\qquad  \left.\hat z\left(\xi+\dfrac{ u(s)}{\mathfrak p^{-1}N}+u(N)\right)\overline{\hat w\left(\xi+\dfrac{ u(s)}{\mathfrak p^{-1}N}+u(N)\right)}\right\}=0.
\end{align*}

\parindent=0mm \vspace{.0in}
{\it{Proof.}} For all $z,w \in l^2(\Lambda)$, the orthogonality of the  systems $\left\{ T_{qN\lambda}z\right\}_{\lambda \in \Lambda}$ and $\left\{ T_{qN\lambda}w\right\}_{\lambda \in \Lambda}$,  is equivalent to
\begin{align*}
0&=\big\langle\, T_{qN\lambda}z, T_{qN\sigma}w\big\rangle\\
&=\left\langle\, \big(T_{qN\lambda}z\big)^{\wedge}, \big(T_{qN\sigma}w\big)^{\wedge}\right\rangle\\
&=\int_\Omega \hat z(\xi)\overline{\hat w(\xi)}\,\overline{\chi_{{\mathfrak p}^{-1}N(\lambda-\sigma)}(\xi)} d\xi\\
&=\int_{{\mathfrak p}{ \mathfrak D}}\left\{\hat z(\xi)\overline{\hat w(\xi)}+\hat z\left(\xi+u(N)\right)\overline{\hat w\left(\xi+u(N)\right)}\right\}\overline{\chi_{{\mathfrak p}^{-1}N(\lambda-\sigma)}(\xi)} d\xi \quad \lambda, \sigma \in \Lambda. \tag{3.4}
\end{align*}

\parindent=0mm \vspace{.0in}
Taking $\lambda = u(m)$ and $\sigma = u(n)$  in (3.4) and setting
\begin{align*}
 h(\xi)=\hat z(\xi)\overline{\hat w(\xi)}+\hat z\left(\xi+u(N)\right)\overline{\hat w\left(\xi+u(N)\right)},\tag{3.5}
 \end{align*}
we obtain
\begin{align*}
0&=\int_{{\mathfrak p}{ \mathfrak D}}h(\xi)\,\overline{\chi\Big(({\mathfrak p}^{-1}N)u\big(q(m-n)\big) \xi\Big)}  d\xi \\
&= \int_{({\mathfrak p}/ qN ){\mathfrak D}}\left\{ \sum_{s=0}^{qN-1}h\left(\xi+\dfrac{ u(s)}{\mathfrak p^{-1}N}\right)\right\}\overline{\chi\Big(({\mathfrak p}^{-1}N)u\big(q(m-n)\big) \xi\Big)}  d\xi.
\end{align*}

\parindent=0mm \vspace{.0in}
Since the above equality holds for all $(m-n)\in \mathbb N_0$, it follows that
\begin{align*}
\sum_{s=0}^{qN-1}h\left(\xi+\dfrac{ u(s)}{\mathfrak p^{-1}N}\right)=0,\quad a.e. \tag{3.6}
\end{align*}

\parindent=0mm \vspace{.0in}
On taking $\lambda = r/N+u(m)$ and $\sigma = r/N+u(n)$ where $m, n \in \mathbb N_{0}$, we obtain the same identity (3.6). Similarly, taking $\lambda = r/N+u(m)$ and $\sigma = u(n)$, where $m, n \in \mathbb N_{0}$, in (3.3), we obtain
\begin{align*}
0&=\int_{{\mathfrak p}{ \mathfrak D}}h(\xi)\,\overline{\chi\Big(({\mathfrak p}^{-1}N)u\big(q(m-n)\big)\xi\Big)}\, \overline{\chi({\mathfrak p}^{-1}r \xi)}\, d\xi \\
&= \int_{({\mathfrak p}/ qN ){\mathfrak D}}\left\{ \sum_{s=0}^{qN-1}\overline{\chi\left(\dfrac{r}{N}\mathfrak p u(s)\right)}h\left(\xi+\dfrac{ u(s)}{\mathfrak p^{-1}N}\right)\right\}\overline{\chi\Big(({\mathfrak p}^{-1}N)u\big(q(m-n)\big)\xi\Big)} \overline{\chi({\mathfrak p}^{-1}r\xi) }\, d\xi.
\end{align*}

\parindent=0mm \vspace{.1in}
Thus, we conclude that
\begin{align*}
\sum_{s=0}^{qN-1}h\left(\xi+\dfrac{ u(s)}{\mathfrak p^{-1}N}\right)\overline{\chi\left(\dfrac{r}{N}\mathfrak p u(s)\right)}=0,\quad a.e. \tag{3.7}
\end{align*}

\parindent=0mm \vspace{.0in}
By taking $\lambda = u(m)$ and $\sigma =\dfrac{r}{N}+u(n)$, we have $\lambda-\sigma=-\dfrac{r}{N}+u(m-n),m, n \in \mathbb N_{0}$ and consequently, we get
\begin{align*}
\sum_{s=0}^{qN-1}h\left(\xi+\dfrac{ u(s)}{\mathfrak p^{-1}N}\right)\chi\left(\dfrac{r}{N}\mathfrak p u(s)\right)=0,\quad a.e. \tag{3.8}
\end{align*}

\parindent=0mm \vspace{.0in}
This completes the proof of the Theorem 3.3.\qquad \fbox

\parindent=8mm \vspace{.1in}

Since the collection $\left\{\chi_\lambda(\xi):\lambda \in \Lambda\right\}$ constitutes an orthonormal basis for $L^2(\Omega)$, there exist locally $L^2, {\mathfrak p}-$periodic functions $z_1, z_2, w_1$ and $w_2$ such that
\begin{align*}
\hat z(\xi)= z_1(\xi)+ \overline{\chi\left(\dfrac{r}{N}\xi\right)}\,z_2(\xi),\quad \text{and} \quad  \hat w(\xi)= w_1(\xi)+ \overline{\chi\left(\dfrac{r}{N}\xi\right)}\,w_2(\xi).\tag{3.9}
\end{align*}

\parindent=0mm \vspace{.0in}
Therefore, equation (3.5) becomes
\begin{align*}
 h(\xi)=q\left(z_1(\xi)\overline{w_1(\xi)}+z_2(\xi)\overline{w_2(\xi)}\right).\tag{3.10}
\end{align*}

\parindent=0mm \vspace{.0in}
Consequently, Eqs.(3.6)--(3.9) yields
\begin{align*}
&\sum_{s=0}^{qN-1}\left\{z_1\left(\xi+\dfrac{ u(s)}{\mathfrak p^{-1}N}\right)\overline{w_1\left(\xi+\dfrac{ u(s)}{\mathfrak p^{-1}N}\right)} z_2\left(\xi+\dfrac{ u(s)}{\mathfrak p^{-1}N}\right)\overline{w_2\left(\xi+\dfrac{ u(s)}{\mathfrak p^{-1}N}\right)}\right\}=0,\tag{3.11}\\
&\sum_{s=0}^{qN-1}\overline{\chi\left(\dfrac{r}{N}\, \mathfrak p u(s) \right)}\left\{z_1\left(\xi+\dfrac{ u(s)}{\mathfrak p^{-1}N}\right)\overline{w_1\left(\xi+\dfrac{ u(s)}{\mathfrak p^{-1}N}\right)} z_2\left(\xi+\dfrac{ u(s)}{\mathfrak p^{-1}N}\right)\overline{w_2\left(\xi+\dfrac{ u(s)}{\mathfrak p^{-1}N}\right)}\right\}=0, \tag{3.12}\\
&\sum_{s=0}^{qN-1}\chi\left(\dfrac{r}{N}\, \mathfrak p u(s) \right)\left\{z_1\left(\xi+\dfrac{ u(s)}{\mathfrak p^{-1}N}\right)\overline{w_1\left(\xi+\dfrac{ u(s)}{\mathfrak p^{-1}N}\right)} z_2\left(\xi+\dfrac{ u(s)}{\mathfrak p^{-1}N}\right)\overline{w_2\left(\xi+\dfrac{ u(s)}{\mathfrak p^{-1}N}\right)}\right\}=0.\tag{3.13}
\end{align*}

\parindent=0mm \vspace{.0in}
{\bf Corollary 3.4.} {\it Suppose for $z,w \in l^2(\Lambda)$ and equation (3.9) is satisfied. Then, the subspaces generated by the systems $\left\{ T_{qN\lambda}z\right\}_{\lambda \in \Lambda}$ and $\left\{ T_{qN\lambda}w\right\}_{\lambda \in \Lambda}$ are orthogonal in $l^2(\Lambda)$ if and only if the equations (3.11)-(3.13) are satisfied.}

\parindent=0mm \vspace{.2in}
{\bf{Theorem 3.5.}} {\it For $z \in l^2(\Lambda)$, the system $\left\{ T_{qN\lambda}z\right\}_{\lambda \in \Lambda}$ is orthonormal  in $l^2(\Lambda)$ if and only if the following conditions hold:}
\begin{align*}
&\dfrac{1}{q} \sum_{s=0}^{qN-1}\left\{\left|\hat z\left(\xi+\dfrac{ u(s)}{\mathfrak p^{-1}N}\right)\right|^2+\left|\hat z\left(\xi+\dfrac{ u(s)}{\mathfrak p^{-1}N}+u(N)\right)\right|^2\right\}=qN,\tag{3.14}\\
&\sum_{s=0}^{qN-1}\overline{\chi\left(\dfrac{r}{N}\, \mathfrak p u(s) \right)}\left\{\left|\hat z\left(\xi+\dfrac{ u(s)}{\mathfrak p^{-1}N}\right)\right|^2+\left|\hat z\left(\xi+\dfrac{ u(s)}{\mathfrak p^{-1}N}+u(N)\right)\right|^2\right\}=0. \tag{3.15}
\end{align*}

\parindent=0mm \vspace{.0in}
{\it{Proof.}} For $z \in l^2(\Lambda)$, and $\lambda, \sigma \in \Lambda$, the orthonormality of the system $\left\{ T_{qN\lambda}z\right\}_{\lambda \in \Lambda}$  in $l^2(\Lambda)$ is equivalent to
\begin{align*}
\int_{{\mathfrak p}{ \mathfrak D}}\Big\{\left|\hat z\left(\xi\right)\right|^2+\left|\hat z\left(\xi+u(N)\right)\right|^2\Big\}\overline{\chi\Big(({\mathfrak p}^{-1}N)u\big(q(m-n)\big)\xi\Big)}d\xi=\delta_{\lambda, \sigma}.
\end{align*}

\parindent=0mm \vspace{.0in}
Proceeding in a similar way as in the proof of Theorem 3.3, we obtain the desired result.\qquad \fbox

\parindent=0mm \vspace{.1in}

{\bf Corollary 3.6.} {\it Let $z \in l^2(\Lambda)$ be such that  condition (3.9) holds. Then, the  system $\left\{ T_{qN\lambda}z\right\}_{\lambda \in \Lambda}$ is orthogonal in $l^2(\Lambda)$ if and only if the  following identities  hold:}
\begin{align*}
&\text{(a)}\quad\sum_{s=0}^{qN-1}\left\{\left|\hat z\left(\xi+\dfrac{ u(s)}{\mathfrak p^{-1}N}\right)\right|^2+\left|\hat z\left(\xi+\dfrac{ u(s)}{\mathfrak p^{-1}N}\right)\right|^2\right\}=qN,\qquad\qquad\qquad\qquad\\
&\text{(b)}\quad \sum_{s=0}^{qN-1}\overline{\chi\left(\dfrac{r}{N}\, \mathfrak p u(s) \right)}\left\{\left|\hat z\left(\xi+\dfrac{ u(s)}{\mathfrak p^{-1}N}\right)\right|^2+\left|\hat z\left(\xi+\dfrac{ u(s)}{\mathfrak p^{-1}N}\right)\right|^2\right\}=0.\qquad\qquad\qquad\qquad
\end{align*}

\parindent=0mm \vspace{.0in}
{\bf{Theorem 3.7.}} {\it If  ${\mathscr F(W)}$ is the first-stage discrete wavelet system as defined by (3.3). Then the following statements are equivalent:}

\parindent=0mm \vspace{.1in}
(i) {\it The set ${\mathscr F(W)}$ is an orthonormal basis for $l^2(\Lambda)$}

\parindent=0mm \vspace{.1in}
(ii) {\it The matrix ${\mathcal M}(\xi)$ of order $q^2N \times q^2N$ is unitary,  when the entries ${\mathcal M}_{st}(\xi),\,0\le s,t \le q^2N-1$ of ${\mathcal M}(\xi)$ are defined as follows:}
\begin{align*}
{\mathcal M}_{st}(\xi)=\dfrac{1}{q\sqrt N}\left\{
\begin{array}{lcr}
\hat w_t\left(\xi+\dfrac{ u(s)}{\mathfrak p^{-1}N}\right);\quad 0\le s \le qN-1;0\le t \le qN-1,\\\\
\hat w_t\left(\xi+\dfrac{u(s-qN)}{\mathfrak p^{-1} N}+u(N)\right);\quad qN\le s \le q^2N-1;\\
\qquad\qquad\qquad \qquad\qquad\qquad\qquad\qquad0\le t \le qN-1,\\\\
\overline{\chi\left(\dfrac{r}{N}\, \mathfrak p u(s) \right)}\hat w_{t-qN}\left(\xi+\dfrac{u(s-qN)}{\mathfrak p^{-1} N}\right);\quad 0\le t \le qN-1;\\
\qquad\qquad\qquad \qquad\qquad\qquad\qquad\quad\qquad qN\le s \le q^2N-1,\\
\overline{\chi\left(\dfrac{r}{N}\, \mathfrak p u(s) \right)}\hat w_{t-qN}\left(\xi+\dfrac{u(s-qN)}{\mathfrak p^{-1} N}+u(N)\right);\\
\qquad\qquad\qquad\quad\qquad\qquad\qquad\qquad qN\le t,s \le q^2N-1.
\end{array}\right.\tag{3.16}
\end{align*}

\parindent=0mm \vspace{.1in}
{\it{Proof.}} Suppose that the system ${\mathscr F(W)}$ defined by (3.3) is an orthonormal basis for $l^2(\Lambda)$. Then, ${\mathscr F(W)}$ is an orthonormal set in $l^2(\Lambda)$ and therefore, for $\lambda, \sigma \in \Lambda$, and $0 \le \ell, k \le qN-1$, we have
\begin{align*}
\big\langle\, T_{qN\lambda}w_\ell, T_{qN\sigma}w_k\big\rangle=\delta_{\lambda, \sigma }\delta_{\ell, k}.\end{align*}

\parindent=0mm \vspace{.0in}
For each $\ell$ and $k,0 \le \ell, k \le qN-1$, Theorems 3.3 and 3.6 implies that
\begin{align*}
 &\dfrac{1}{q}\displaystyle\sum_{t=0}^{qN-1}\left\{\hat w_\ell \left(\xi+\dfrac{u(t)}{\mathfrak p^{-1}N}\right)\overline{\hat w_k\left(\xi+\dfrac{u(t)}{\mathfrak p^{-1}N}\right)}\right.\\
&\qquad\qquad+\left.\hat w_\ell\left(\xi+\dfrac{u(t)}{\mathfrak p^{-1}N}+u(N)\right)\overline{\hat w_k\left(\xi+\dfrac{u(t)}{\mathfrak p^{-1}N}+u(N)\right)}\right\}=qN \delta_{\ell, k},\tag{3.17}\\
&{\mathcal R}_{\ell, k}(\xi)=\sum_{t=0}^{qN-1}\overline{\chi\left(\dfrac{r}{N}\, \mathfrak p u(t) \right)}\left\{\hat w_\ell\left(\xi+\dfrac{u(t)}{\mathfrak p^{-1}N}\right)\overline{\hat w_k\left(\xi+\dfrac{u(t)}{\mathfrak p^{-1}N}\right)}\right.\\
&\qquad\qquad\qquad+\left.\hat w_\ell\left(\xi+\dfrac{u(t)}{\mathfrak p^{-1}N}+u(N)\right)\overline{\hat w_k\left(\xi+\dfrac{u(t)}{\mathfrak p^{-1}N}+u(N)\right)}\right\}=0, \tag{3.18}
\end{align*}
and
\begin{align*}
 & {\mathcal S}_{\ell, k}(\xi)=\displaystyle\sum_{t=0}^{qN-1}\chi\left(\dfrac{r}{N}\, \mathfrak p u(t) \right)\left\{\hat w_\ell\left(\xi+\dfrac{u(t)}{\mathfrak p^{-1}N}\right)\overline{\hat w_k\left(\xi+\dfrac{u(t)}{\mathfrak p^{-1}N}\right)}\right.\\
&\qquad\qquad+\left.\hat w_\ell\left(\xi+\dfrac{u(t)}{\mathfrak p^{-1}N}+u(N)\right)\overline{\hat w_k\left(\xi+\dfrac{u(t)}{\mathfrak p^{-1}N}+u(N)\right)}\right\}=0. \tag{3.19}
\end{align*}

\parindent=0mm \vspace{.0in}
Thus, it is sufficient to consider the equations (3.17) and (3.18), as ${\mathcal S}_{\ell, k}(\xi)=\overline{{\mathcal R}_{\ell, k}(\xi)}$. These equations give rise to the matrix ${\mathcal M}(\xi)$ of order $q^2N \times q^2N$ with entries as defined in system (3.164). Moreover, the matrix ${\mathcal M}(\xi)$ is unitary. This follows from the identities (3.17) and (3.18) and noting that columns of ${\mathcal M}(\xi)$ form an orthonormal system in ${\mathbb C}^{q^2N}$ with respect to the usual inner product, and hence form an orthonormal basis of ${\mathbb C}^{q^2N}$.

\parindent=8mm \vspace{.1in}
Conversely, assume that the matrix ${\mathcal M}(\xi)$ is unitary. Then, it suffices to show that the set ${\mathscr F(W)}$ is complete. For this, let $w \in l^2(\Lambda)$ and $Pw$ be the projection onto $\overline{\text {span}}\,{\mathscr F(W)}$, then we have
\begin{align*}
\big\|Pw\big\|^2_2&=\sum_{\ell=0}^{qN-1}\sum_{\lambda \in \Lambda}\left|\big\langle w, T_{qN\lambda}w_\ell\big\rangle\right|^2\\
&=\sum_{\ell=0}^{qN-1}\sum_{\lambda \in \Lambda}\left|\big\langle \hat w, (T_{qN\lambda}w_\ell)^{\wedge}\big\rangle\right|^2\\
&=\sum_{\ell=0}^{qN-1}\sum_{\lambda \in \Lambda}\left|\int_{\Omega }\hat w(\xi)\overline{\hat w_\ell(\xi)}\chi_\lambda({\mathfrak p}^{-1}N\xi)d\xi\right|^2.
\end{align*}

\parindent=0mm \vspace{.0in}
Writing $\Omega={\mathfrak p}{\mathfrak D}\cup{\mathfrak p}(N+{\mathfrak D})$ and observing that ${\mathfrak p}{\mathfrak D}=\bigcup_{s=0}^{qN-1}{\mathfrak p}(s+{\mathfrak D})/N$ , we obtain

\parindent=0mm \vspace{.1in}
$\big\|Pw\big\|^2_2$
\begin{align*}
=\dfrac{1}{(qN)^2}\displaystyle\sum_{\ell=0}^{qN-1}\sum_{\lambda \in \Lambda}\left|\int_{{\mathfrak p}{\mathfrak D}}\sum_{t=0}^{qN-1}\chi_\lambda\big( \mathfrak p u(t) \big)\left\{\hat w_\ell\left(\zeta_t\right)\overline{\hat w_\ell\left(\zeta_t\right)}+\hat w\left(\zeta_t+u(N)\right)\overline{\hat w_\ell\left(\zeta_t\right)}\right\}\chi_\lambda(\xi)d\xi\right|^2.\tag{3.20}
\end{align*}
where $\zeta_t=\big(\xi+ u(t)\big)/\mathfrak p^{-1}N$. As $\Lambda = {\mathcal Z} \cup \left\{ {\mathcal Z}+ r/N\right\}$, we can rewrite (3.20), by using Plancherel formula as follows

\parindent=0mm \vspace{.1in}
$\big\|Pw\big\|^2_2$
\begin{align*}
&=\dfrac{1}{(qN)^2}\displaystyle\sum_{\ell=0}^{qN-1}\sum_{m \in \mathbb N_0}\left|\int_{{\mathfrak p}{\mathfrak D}}\sum_{s=0}^{qN-1}\left\{\hat w\left(\zeta_s\right)\overline{\hat w_\ell\left(\zeta_s\right)}+\hat w\left(\zeta_s\right)\overline{\hat w_\ell\big(\zeta_s+u(N)\big)}\right\}\chi\big({\mathfrak p}u(m)\xi\big)d\xi\,\right|^2\\
&\qquad\qquad\qquad+\dfrac{1}{(qN)^2}\displaystyle\sum_{\ell=0}^{qN-1}\sum_{m \in \mathbb N_0}\left|\int_{{\mathfrak p}{\mathfrak D}}\sum_{s=0}^{qN-1}\chi\left(\dfrac{r}{N}\, \mathfrak p u(s) \right)\left\{\hat w\left(\zeta_s\right)\overline{\hat w_\ell\left(\zeta_s\right)}\right.\right.\\
&\left.\left.\qquad\qquad\qquad+\hat w\left(\zeta_s+u(N)\right)\overline{\hat w_\ell\big(\zeta_s+u(N)\big)}\chi\left(\dfrac{r}{N}\, \xi \right)\right\}\chi\big({\mathfrak p}u(m)\xi\big)d\xi\,\right|^2\\
&=\dfrac{1}{q^3N^2}\displaystyle\sum_{\ell=0}^{qN-1}\left|\int_{{\mathfrak p}{\mathfrak D}}\sum_{s=0}^{qN-1}\left\{\hat w\left(\zeta_s\right)\overline{\hat w_\ell\left(\xi+\dfrac{ u(s)}{\mathfrak p^{-1}N}\right)}\right.+\hat w\big(\zeta_s+u(N)\big)\overline{\hat w_\ell\big(\zeta_s+u(N)\big)}\Big\}\right|^2d\xi\\
&+\dfrac{1}{q^3N^2}\displaystyle\sum_{\ell=0}^{qN-1}\left|\int_{{\mathfrak p}{\mathfrak D}}\sum_{s=0}^{qN-1}\chi\left(\dfrac{r}{N}\, \mathfrak p u(s) \right)\left\{\hat w\left(\zeta_s\right)\overline{\hat w_\ell\left(\zeta_s\right)}\right.+\hat w\big(\zeta_s+u(N)\big)\overline{\hat w_\ell\big(\zeta_s+u(N)\big)}\right\}\Big|^2d\xi.
\end{align*}

\parindent=0mm \vspace{.0in}
Since the rows of matrix ${\mathcal M}(\xi)$ form an orthonormal system in ${\mathbb C}^{q^2N}$, therefore, for each $0 \le r,s \le qN-1$, we have \begin{align*}
&\dfrac{1}{q}\sum_{\ell=0}^{qN-1}\left\{\left(1-\overline{\chi\left(\dfrac{r}{N}\, \mathfrak p u(r)- \mathfrak p u(s)\right)}\right)\hat w_\ell\left(\xi+\dfrac{ u(s)}{\mathfrak p^{-1}N}\right)\overline{\hat w_\ell\left(\xi+\dfrac{ u(s)}{\mathfrak p^{-1}N}\right)}\right\}=qN \delta_{r,s}\\
&\dfrac{1}{q}\sum_{\ell=0}^{qN-1}\left\{\left(1-\overline{\chi\left(\dfrac{r}{N}\, \mathfrak p u(r)- \mathfrak p u(s)\right)}\right)\hat w_\ell\left(\xi+\dfrac{ u(s)}{\mathfrak p^{-1}N}+u(N)\right)\overline{\hat w_\ell\left(\xi+\dfrac{ u(s)}{\mathfrak p^{-1}N}+u(N)\right)}\right\}=qN \delta_{r,s}
\end{align*}
and
\begin{align*}
\sum_{\ell=0}^{qN-1}\left\{\left(1-\overline{\chi\left(\dfrac{r}{N}\, \mathfrak p u(r)- \mathfrak p u(s)\right)}\right)\hat w_\ell\left(\xi+\dfrac{ u(s)}{\mathfrak p^{-1}N}\right)\overline{\hat w_\ell\left(\xi+\dfrac{ u(s)}{\mathfrak p^{-1}N}+u(N)\right)}\right\}=0.
\end{align*}

\parindent=0mm \vspace{.1in}
Thus, for $0 \le s,t \le qN-1$, we have
\begin{align*}
\big\|Pw\big\|^2_2&=\dfrac{1}{(qN)^2}\int_{{\mathfrak p}{\mathfrak D}}\sum_{\ell=0}^{qN-1}\sum_{s=0}^{qN-1}\left\{\left|\hat w\left(\zeta_s\right)\overline{\hat w_\ell\left(\zeta_s\right)}\right|^2+\left|\hat w\left(\zeta_s+u(N)\right)\overline{\hat w_\ell\left(\zeta_s+u(N)\right)}\right|^2\right\}d\xi\\
&=\dfrac{1}{(qN)^2}\int_{{\mathfrak p}{\mathfrak D}}\sum_{s=0}^{qN-1}\left\{\left|\hat w\left(\zeta_s\right)\right|^2\sum_{\ell=0}^{qN-1}\left|\overline{\hat w_\ell\left(\zeta_s\right)}\right|^2+\Big|\hat w\big(\zeta_s+u(N)\big)\Big|^2\displaystyle\sum_{\ell=0}^{qN-1}\left|\overline{\hat w_\ell\left(\zeta_s+u(N)\right)}\right|^2\right\}d\xi\\
&=\dfrac{1}{qN}\int_{{\mathfrak p}{\mathfrak D}}\sum_{s=0}^{qN-1}\left\{\left|\hat w\left(\zeta_s\right)\right|^2+\left|\hat w\big(\zeta_s+u(N)\big)\right|^2\right\}d\xi\\
&=\dfrac{1}{qN}\sum_{s=0}^{qN-1}\int_{(1+s){\mathfrak p}{\mathfrak D}}\left\{\left|\hat w\left(\dfrac{\mathfrak p \xi }{N}\right)\right|^2+\left|\hat w\left(\dfrac{\mathfrak p \xi }{N}+u(N)\right)\right|^2\right\}d\xi\\
&=\dfrac{1}{qN}\int_{N{\mathfrak D}}\left\{\left|\hat w\left(\dfrac{\mathfrak p \xi }{N}\right)\right|^2+\left|\hat w\left(\dfrac{\mathfrak p \xi }{N}+u(N)\right)\right|^2\right\}d\xi\\
&=\int_{{\mathfrak p}{\mathfrak D}}\left\{\left|\hat w\left(\xi\right)\right|^2+\left|\hat w\big(\xi+u(N)\big)\right|^2\right\}d\xi\\
&=\int_{\Omega}\left|\hat w\left(\xi\right)\right|^2d\xi\\
&=\big\| w\big\|^2_2.
\end{align*}

\parindent=0mm \vspace{.0in}
Hence, the projection $P$ is an identity map and $\overline{\text {span}}\,{\mathscr F(W)}=l^2(\Lambda)$. Therefore, the set $\mathcal F$ is an orthonormal basis for $l^2(\Lambda)$. This completes the proof of Theorem 3.7. \qquad \fbox

\parindent=0mm \vspace{.1in}

{\bf{Corollary 3.8.}} {\it For each $\ell,0 \le \ell \le qN-1$, let $w_\ell \in l^2(\Lambda)$ such that
\begin{align*}
\hat w_\ell(\xi)= w_{\ell_0}(\xi)+ \overline{\chi\left(\dfrac{r}{N}\xi\right)}\,w_{\ell_1}(\xi),\tag{3.21}
\end{align*}

\parindent=0mm \vspace{.0in}
for some ${\mathfrak p}$-periodic functions $ w_{\ell_0}$ and $w_{\ell_1}$. Then, the system ${\mathscr F(W)}$ as defined by (3.3) is an orthonormal basis for $l^2(\Lambda)$ if and only if the matrix
\begin{align*}
{\mathcal M}_{st}(\xi)=\dfrac{1}{\sqrt{q N}}\left\{
\begin{array}{lcr}
 w_{t_0}\left(\xi+\dfrac{ u(s)}{\mathfrak p^{-1}N}\right);\quad 0\le s \le qN-1;0\le t \le qN-1,&&\\
 w_{t_1}\left(\xi+\dfrac{u(s-qN)}{\mathfrak p^{-1} N}\right);\quad qN\le s \le q^2N-1;0\le t \le qN-1,&&\\
\overline{\chi\left(\dfrac{r}{N}\, \mathfrak p u(s) \right)} w_{(t-qN)0}\left(\xi+\dfrac{u(s-qN)}{\mathfrak p^{-1} N}\right);\quad 0\le t \le qN-1;&&\\
\qquad\qquad\qquad \qquad\qquad\qquad\qquad\qquad\qquad qN\le s \le q^2N-1,&&\\
\overline{\chi\left(\dfrac{r}{N}\, \mathfrak p u(s) \right)} w_{(t-qN)1}\left(\xi+\dfrac{u(s-qN)}{\mathfrak p^{-1} N}\right);\quad qN\le t,s \le q^2N-1;
\end{array}\right.
\end{align*}
is unitary.}

\parindent=8mm \vspace{.0in}
Assume $w_0$ in $l^2(\Lambda)$ such that it satisfies equations (3.14) and (3.15). Following the procedure of the paper of Shah and Abdullah \cite{ShahA3}, it can be easily  shown that there exists functions $w_k: 1\le k \le qN-1$ satisfying conditions (3.17) and (3.18) if and only if the function $M_0$ is of the form
\begin{align*}
M_0(\xi)=\left|\frac{\hat w_0(\xi)}{q\sqrt N}\right|^2+\left|\frac{\hat w_0\big(\xi+u(N)\big)}{q\sqrt N}\right|^2,
\end{align*}
and satisfies the following identity
\begin{align*}
M_0\big(\xi+{\mathfrak p}^2\big)=M_0(\xi).
\end{align*}

\parindent=0mm \vspace{.1in}
{\bf{Theorem 3.9.}} {\it For each $i \in \{0,1,\dots,qN-1\}$ and $\ell \in \mathbb N$, let $f_{\ell,i}\in l^2(\Lambda)$. Then, the system
\begin{align*}
{\mathcal G}(H)=\Big\{T_{(qN)^\ell\lambda}h_{\ell,i}: \lambda \in \Lambda , 0 \le i \le qN-1\Big\}\tag{3.22}
\end{align*}

\parindent=0mm \vspace{.0in}
is an orthonormal in $l^2(\Lambda)$ if and only if for $0 \le i,j \le qN-1$, the following conditions are satisfied:}
\begin{align*}
 &{\mathcal T}_{ij}^\ell(\xi)=\displaystyle\sum_{s=0}^{(qN)^\ell-1}\left\{\hat h_{\ell,i}\left(\xi+\dfrac{u(s)}{(\mathfrak p^{-1}N)^\ell}\right)\overline{\hat h_{\ell,j}\left(\xi+\dfrac{u(s)}{(\mathfrak p^{-1}N)^\ell}\right)}\right.\\
&\qquad+\left.\hat h_{\ell,i}\left(\xi+\dfrac{u(s)}{(\mathfrak p^{-1}N)^\ell}+u(N)\right)\overline{\hat h_{\ell,j}\left(\xi+\dfrac{u(s)}{(\mathfrak p^{-1}N)^\ell}+u(N)\right)}\right\}=q(qN)^\ell\delta_{i,j},\tag{3.23}\\
&\sum_{s=0}^{(qN)^\ell-1}\overline{\chi\left(\dfrac{r}{N}\, \mathfrak p u(s) \right)}\left\{\hat h_{\ell,i}\left(\xi+\dfrac{u(s)}{(\mathfrak p^{-1}N)^\ell}\right)\overline{\hat h_{\ell,j}\left(\xi+\dfrac{u(s)}{(\mathfrak p^{-1}N)^\ell}\right)}\right.\\
&\qquad+\left.\hat h_{\ell,i}\left(\xi+\dfrac{u(s)}{(\mathfrak p^{-1}N)^\ell}+u(N)\right)\overline{\hat h_{\ell,j}\left(\xi+\dfrac{u(s)}{(\mathfrak p^{-1}N)^\ell}+u(N)\right)}\right\}=0.\tag{3.24}
\end{align*}

\parindent=0mm \vspace{.1in}
{\it{Proof.}} For $\lambda, \sigma \in \Lambda$ and $0 \le i,j \le qN-1$, the orthonormality of the system ${\mathcal G}(H)$ in $l^2(\Lambda)$ is equivalent to
\begin{align*}
\big\langle\, T_{(qN)^\ell\lambda}h_{\ell,i}, T_{(qN)^\ell\sigma}h_{\ell,j}\big\rangle=\delta_{\lambda, \sigma }\delta_{i,j}.
\end{align*}

\parindent=0mm \vspace{.0in}
By setting
\begin{align*}  H_{ij}(\xi)=\hat h_{\ell,i}(\xi)\overline{\hat h_{\ell,j}(\xi)}+\hat h_{\ell,i}\big(\xi+u(N)\big)\overline{\hat h_{\ell,j}\big(\xi+u(N)\big)},
\end{align*}

\parindent=0mm \vspace{.0in}
and taking $\lambda=u(m)$ and $\sigma=u(n)$, where $m,n \in \mathbb N_0$, we have
\begin{align*}
\delta_{\lambda, \sigma }\delta_{i,j}&=\int_{{\mathfrak p}{\mathfrak D}}H_{ij}(\xi)\overline{\chi\Big({\mathfrak p}^{-1}({\mathfrak p}^{-1}N)^\ell\big(u(m)-u(n)\big) \xi\Big)}  d\xi \\
&= \int_{({\mathfrak p}/( qN)^\ell ){\mathfrak D}}\left\{ \sum_{s=0}^{(qN)^\ell-1}H_{ij}\left(\xi+\dfrac{u(s)}{(\mathfrak p^{-1}N)^\ell}\right)\right\}\overline{\chi\Big(({\mathfrak p}^{-1}N)^\ell u\big(q(m-n)\big) \xi\Big)}  d\xi.
\end{align*}

\parindent=0mm \vspace{.0in}
Now the desired result can be proved analogously to Theorem  3.3. \qquad \fbox

\parindent=0mm \vspace{.2in}

{\bf{4. $J^{\text{th}}$-stage Discrete Wavelets on local Fields}}

\parindent=0mm \vspace{.1in}
In this section, we introduce the notion of $J^{\text{th}}$-stage nonuniform discrete wavelet system in the Hilbert space $l^2(\Lambda)$ and show that this space can be expressed as an orthogonal decomposition in terms of countable number of its closed subspaces.

\parindent=0mm \vspace{.1in}
{\bf{Definition 4.1.}}   Let $N \in\mathbb N$ and let $k$ be an odd integer with $1 \le k \le qN-1$ such that $k$ and $N$ are relatively prime. For $h_{j,k}\in l^2(\Lambda), J\in \mathbb N$, we call $\mathscr H$ a $J^{\text{th}}$-stage nonuniform discrete wavelet system associated with $H=\left\{ h_{j,k}: h_{j,k}\in l^2(\Lambda)\right\}$ if
\begin{align*}
{\mathscr H}=\Big\{T_{(qN)^{j}\lambda}h_{j,k}: h_{j,k}\in l^2(\Lambda); 1\le j\le J,  1 \le k \le qN-1,\lambda \in \Lambda\Big\}. \tag{4.1}
\end{align*}
 is a complete orthonormal set in $l^2(\Lambda)$.

\parindent=0mm \vspace{.1in}
{\bf{Theorem 4.2.}} {\it Let the system $\big\{T_{qN\lambda}w_i: w_i\in l^2(\Lambda); \lambda \in \Lambda ; 0 \le i \le qN-1\big\}$ be orthonormal in $l^2(\Lambda )$, where $w_i\in l^2(\Lambda)$ satisfies equation (3.21). For $\ell \in \mathbb N$ and $h_{(\ell-1),i}\in l^2(\Lambda)$, let the system
$\big\{T_{qN\lambda}h_{(\ell-1),i}: \lambda \in \Lambda , 0 \le i \le qN-1\big\}$ be orthonormal in $l^2(\Lambda )$. Consider the following relation
\begin{align*}
\hat h_{\ell,i}(\xi)=\hat h_{(\ell-1),i}(\xi)\hat w_i\big(({\mathfrak p}^{-1}N)^{\ell-1}\xi\big),\tag{4.2}
\end{align*}
where $\hat h_{0,0}(\xi)=1$ a.e. Then, the system
\begin{align*}
\Big\{T_{(qN)^\ell\lambda}h_{\ell,i}:  \lambda \in \Lambda ; 0 \le i \le qN-1\Big\}\tag{4.3}
\end{align*}
is orthonormal in $l^2(\Lambda ).$}

\parindent=0mm \vspace{.1in}
{\it{Proof.}} For $\ell=1,$ the result follows immediately. To prove the required result for $\ell \in \mathbb N -\{1\}$ and $0 \le i, j \le qN-1$, it is sufficient to prove the identities (3.23) and (3.24). However, we observe that $\hat h_{\ell,i} \in L^2(\Omega)$ since
\begin{align*}
\left\| \hat h_{\ell,i}\right\|^2_2&=\int_\Omega\left|\hat h_{(\ell-1),i}(\xi)\hat w_i\big(({\mathfrak p}^{-1}N)^{\ell-1}\xi\big)\right|^2d\xi\\
&\le {\sup}_\xi\Big|\hat w_i\big(({\mathfrak p}^{-1}N)^{\ell-1}\xi\big)\Big|^2\left\| \hat h_{(\ell-1),i}\right\|^2_2\\
&= qN\left\| \hat h_{(\ell-1),i}\right\|^2_2.
\end{align*}
From (3.23), we infer that
\begin{align*}
 {\mathcal T}_{ij}^\ell(\xi)&=\sum_{s=0}^{(qN)^\ell-1}\left\{\hat h_{\ell,i}\left(\xi+\dfrac{u(s)}{(\mathfrak p^{-1}N)^\ell}\right)\overline{\hat h_{\ell,j}\left(\xi+\dfrac{u(s)}{(\mathfrak p^{-1}N)^\ell}\right)}\right.&&\\
&\qquad\qquad\qquad+\left.\hat h_{\ell,i}\left(\xi+\dfrac{u(s)}{(\mathfrak p^{-1}N)^\ell}+u(N)\right)\overline{\hat h_{\ell,j}\left(\xi+\dfrac{u(s)}{(\mathfrak p^{-1}N)^\ell}+u(N)\right)}\right\}&&\\
 &=\sum_{m=0}^{qN-1}\sum_{n=0}^{(qN)^\ell-1}\left\{\left|\hat h_{(\ell-1),0}\left(\xi+\dfrac{u(n)}{(\mathfrak p^{-1}N)^{\ell-1}}+\dfrac{u(m)}{(\mathfrak p^{-1}N)^\ell}\right)\right|^2\right.&&\\
&\qquad\qquad\qquad+\left.\left|\hat h_{(\ell-1),0}\left(\xi+\dfrac{u(n)}{(\mathfrak p^{-1}N)^{\ell-1}}+\dfrac{u(m)}{(\mathfrak p^{-1}N)^\ell}+u(N)\right)\right|^2\right\}&&
\end{align*}
\begin{align*}
&\qquad\times\left\{\hat w_i\left(({\mathfrak p}^{-1}N)^{\ell-1}\xi+{\mathfrak p}u(n)+\dfrac{u(m)}{\mathfrak p^{-1}N}\right)\overline{\hat w_i\left(({\mathfrak p}^{-1}N)^{\ell-1}\xi+{\mathfrak p}u(n)+\dfrac{u(m)}{\mathfrak p^{-1}N}\right)}\right\},
\end{align*}

\parindent=0mm \vspace{.0in}
where we have used the fact $\hat w_i\big(\xi+u(N)\big)=\hat w_i(\xi)$ and  (4.3).

\parindent=8mm \vspace{.1in}
 For $0\le m \le qN-1$ and $0\le n \le (qN)^{\ell-1}-1$, we define
\begin{align*}
&{H}_{m,n}(\xi)\\
&=\left|\hat h_{(\ell-1),0}\left(\xi+\dfrac{u(n)}{(\mathfrak p^{-1}N)^{\ell-1}}+\dfrac{u(m)}{(\mathfrak p^{-1}N)^\ell}\right)\right|^2+\left|\hat h_{(\ell-1),0}\left(\xi+\dfrac{u(n)}{(\mathfrak p^{-1}N)^{\ell-1}}+\dfrac{u(m)}{(\mathfrak p^{-1}N)^\ell}+u(N)\right)\right|^2.
\end{align*}

\parindent=0mm \vspace{.0in}
Then, we can write
\begin{align*}
 {\mathcal T}_{ij}^\ell(\xi)&=\sum_{m=0}^{qN-1}\sum_{n=0}^{(qN)^\ell-1}{H}_{m,n}(\xi)\left\{\hat w_i\left(({\mathfrak p}^{-1}N)^{\ell-1}\xi+{\mathfrak p}u(n)+\dfrac{u(m)}{\mathfrak p^{-1}N}\right)\right.\\
 &\qquad\qquad\quad\qquad\qquad\qquad\qquad\quad\left.\times\overline{\hat w_i\left(({\mathfrak p}^{-1}N)^{\ell-1}\xi+{\mathfrak p}u(n)+\dfrac{u(m)}{\mathfrak p^{-1}N}\right)}\right\}.
\end{align*}

\parindent=0mm \vspace{.0in}
Using (3.21) for each $0 \le i, j \le qN-1$, we obtain
\begin{align*}
 {\mathcal T}_{ij}^\ell(\xi)&=\sum_{m=0}^{qN-1}\hat w_{i0}\left(({\mathfrak p}^{-1}N)^{\ell-1}\xi+\dfrac{u(m)}{\mathfrak p^{-1}N}\right)\overline{\hat w_{i0}\left(({\mathfrak p}^{-1}N)^{\ell-1}\xi+\dfrac{u(m)}{\mathfrak p^{-1}N}\right)}\sum_{n=0}^{(qN)^\ell-1}{H}_{m,n}(\xi)\\
 &\;+\sum_{m=0}^{qN-1}\hat w_{i1}\left(({\mathfrak p}^{-1}N)^{\ell-1}\xi+\dfrac{u(m)}{\mathfrak p^{-1}N}\right)\overline{\hat w_{i1}\left(({\mathfrak p}^{-1}N)^{\ell-1}\xi+\dfrac{u(m)}{\mathfrak p^{-1}N}\right)}\sum_{n=0}^{(qN)^\ell-1}{H}_{m,n}(\xi)\\
&\;+\sum_{m=0}^{qN-1}\left\{\overline{\chi\left(\dfrac{r}{N}\left(({\mathfrak p}^{-1}N)^{\ell-1}+\dfrac{u(m)}{\mathfrak p^{-1}N}\right)\right)}\hat w_{i1}\left(({\mathfrak p}^{-1}N)^{\ell-1}\xi+\dfrac{u(m)}{\mathfrak p^{-1}N}\right)\right.\\
&\;\left.\times\overline{\hat w_{j0}\left(({\mathfrak p}^{-1}N)^{\ell-1}\xi+\dfrac{u(m)}{\mathfrak p^{-1}N}\right)}\right\}\left\{\displaystyle\sum_{n=0}^{(qN)^\ell-1}\left(\overline{\chi\left(\dfrac{r}{N}{\mathfrak p}u(n)\right)}\right){H}_{m,n}(\xi)\right\}\\
&\;+\sum_{m=0}^{qN-1}\left\{\overline{\chi\left(\dfrac{r}{N}\left(({\mathfrak p}^{-1}N)^{\ell-1}+\dfrac{u(m)}{\mathfrak p^{-1}N}\right)\right)}\hat w_{i0}\left(({\mathfrak p}^{-1}N)^{\ell-1}\xi+\dfrac{u(m)}{\mathfrak p^{-1}N}\right)\right.\\
&\;\left.\times\overline{\hat w_{j1}\left(({\mathfrak p}^{-1}N)^{\ell-1}\xi+\dfrac{u(m)}{\mathfrak p^{-1}N}\right)}\right\}\left\{\displaystyle\sum_{n=0}^{(qN)^\ell-1}\left(\overline{\chi\left(\dfrac{r}{N}{\mathfrak p}u(n)\right)}\right){H}_{m,n}(\xi)\right\}.
\end{align*}

\parindent=0mm \vspace{.0in}
Theorem 3.9 and the orthonormality property of the system (4.2) further yields
\begin{align*}
\sum_{n=0}^{(qN)^\ell-1}{H}_{m,n}(\xi)=q(qN)^{\ell-1}\quad \text{and} \quad \sum_{n=0}^{(qN)^\ell-1}\left(\overline{\chi\left(\dfrac{r}{N}{\mathfrak p}u(n)\right)}\right){H}_{m,n}(\xi)=0,
\end{align*}
which in turn implies
\begin{align*}
 {\mathcal T}_{ij}^\ell(\xi)&=q(qN)^{\ell-1}\sum_{m=0}^{qN-1}\left\{\hat w_{i0}\left(({\mathfrak p}^{-1}N)^{\ell-1}\xi+\dfrac{u(m)}{\mathfrak p^{-1}N}\right)\overline{\hat w_{i0}\left(({\mathfrak p}^{-1}N)^{\ell-1}\xi+\dfrac{u(m)}{\mathfrak p^{-1}N}\right)}\right.\\
 &\qquad \qquad \qquad\left.+~\hat w_{i1}\left(({\mathfrak p}^{-1}N)^{\ell-1}\xi+\dfrac{u(m)}{\mathfrak p^{-1}N}\right)\overline{\hat w_{i1}\left(({\mathfrak p}^{-1}N)^{\ell-1}\xi+\dfrac{u(m)}{\mathfrak p^{-1}N}\right)}\right\}.
\end{align*}

\parindent=0mm \vspace{.0in}
Note that ${\mathcal T}_{ij}^\ell(\xi)=q(qN)^{\ell-1}(qN \delta_{ij})=q(qN)^\ell \delta_{ij},$ as the system ${\mathscr F}(W)$ given by (3.3) is orthonormal in $l^2(\Lambda )$. This proves the equation (3.23). Similarly, we can prove  (3.24). This completes the proof of the Theorem 4.2. \qquad\fbox

\parindent=8mm \vspace{.1in}

We now invoke Theorem 4.2 to prove the orthogonal splitting properties of the subspaces  $V_j$'s.

\parindent=0mm \vspace{.1in}
{\bf{Theorem 4.3.}} {\it With the assumptions of Theorem 4.2, let us define the subsets $V_{\ell-1}, V_\ell$ and $W_\ell$ of $V_0=l^2(\Lambda )$ by}
\begin{align*}
V_{\ell-1}&=\overline{\text{span}}\big\{T_{(qN)^{\ell-1}\lambda}h_{(\ell-1),0}: \lambda \in \Lambda\big\},\\
V_\ell&=\overline{\text{span}}\big\{T_{(qN)^\ell\lambda}h_{\ell,0}: \lambda \in \Lambda\big\},\\
W_\ell&=\overline{\text{span}}\big\{T_{(qN)^\ell\lambda}h_{\ell,i}: \lambda \in \Lambda, 1 \le i \le qN-1\big\}.
\end{align*}

\parindent=0mm \vspace{.0in}
{\it Then, $V_\ell \oplus W_\ell=V_{\ell-1}$.}

\parindent=0mm \vspace{.1in}
{\it{Proof.}} For each $ i =0,1,\dots, qN-1$ and $\ell \in \mathbb N$, we can write
\begin{align*}
\hat w_i\big(({\mathfrak p}^{-1}N)^{\ell-1}\xi\big)&=\sum_{\nu \in \Lambda}w_i(\nu)\overline{\chi\big(({\mathfrak p}^{-1}N)^{\ell-1}\nu\xi\big)}\\
&=\sum_{\sigma \in \Lambda}w_i\big(\sigma-{\mathfrak p}^{-1}N \lambda\big)\overline{\chi\left(({\mathfrak p}^{-1}N)^{\ell-1}\big(\sigma-{\mathfrak p}^{-1}N \lambda\big)\xi\right)}.
\end{align*}
Therefore, we have
\begin{align*}
&\overline{\chi\big(({\mathfrak p}^{-1}N)^\ell\lambda\xi\big)}\hat w_i\big(({\mathfrak p}^{-1}N)^{\ell-1}\xi\big)\hat h_{(\ell-1),0}(\xi)=\sum_{\sigma \in \Lambda}w_i\big(\sigma-{\mathfrak p}^{-1}N \lambda\big)\overline{\chi\left(({\mathfrak p}^{-1}N)^{\ell-1}\sigma\xi\right)}\hat h_{(\ell-1),0}(\xi),\\
&\text{or}~~\left(T_{(qN)^\ell\lambda}h_{\ell,i}\right)^{\wedge}(\xi)=\sum_{\sigma \in \Lambda}w_i\big(\sigma-{\mathfrak p}^{-1}N \lambda\big)\left(T_{(qN)^{\ell-1}\sigma}h_{\ell-1,i}\right)^{\wedge}(\xi),\\
&\text{or}~~T_{(qN)^\ell\lambda}h_{\ell,i}(\xi)=\sum_{\sigma \in \Lambda}w_i\big(\sigma-{\mathfrak p}^{-1}N \lambda\big)T_{(qN)^{\ell-1}\sigma}h_{\ell-1,i}(\xi),
\end{align*}
which implies $V_\ell$ and $W_\ell$ are the subspaces of $V_{\ell-1}$. Using the facts that: $V_\ell$ is orthogonal to $W_\ell$; $\big\{T_{(qN)^\ell\lambda}h_{\ell,0}: \lambda \in \Lambda\big\}$ and $\big\{T_{(qN)^\ell\lambda}h_{\ell,i}: \lambda \in \Lambda, 1 \le i \le qN-1\big\}$ are orthogonal to each other; and $V_\ell \oplus W_\ell \subset V_{\ell-1}, \ell\in\mathbb N$,  it only needs  to show that $V_{\ell-1}\subset V_\ell \oplus W_\ell $. Thus, we have
\begin{align*}
&T_{(qN)^{\ell-1}\lambda}h_{(\ell-1),0}(\xi)\qquad\qquad\\\
&=\sum_{\sigma \in \Lambda}w_i\big(\sigma-{\mathfrak p}^{-1}N \lambda\big)T_{(qN)^\ell\sigma}f_{\ell-1,0}(\xi)\\
&=\sum_{\sigma \in \Lambda}\left\{\sum_{i=0}^{qN-1}\sum_{\nu \in \Lambda}\Big\langle w_i\big(\sigma-{\mathfrak p}^{-1}N \lambda\big), T_{qN\nu}w_i \Big\rangle T_{qN\nu}w_i(\sigma)\right\}T_{(qN)^\ell\sigma}h_{\ell-1,0}(\xi)\\
&=\sum_{i=0}^{qN-1}\sum_{\nu \in \Lambda}\Big\langle w_i\big(\sigma-{\mathfrak p}^{-1}N \lambda\big), T_{qN\nu}w_i \Big\rangle\left\{\sum_{\sigma \in \Lambda} T_{qN\nu}w_i(\sigma)T_{(qN)^\ell\sigma}h_{\ell-1,0}(\xi)\right\}\\
&=\sum_{i=0}^{qN-1}\sum_{\nu \in \Lambda}\Big\langle w_i\big(\sigma-{\mathfrak p}^{-1}N \lambda\big), T_{qN\nu}w_i \Big\rangle T_{(qN)^\ell\nu}h_{\ell, i}(\xi)\\
&=\sum_{\nu \in \Lambda}\Big\langle w_i\big(\sigma-{\mathfrak p}^{-1}N \lambda\big), T_{qN\nu}w_i \Big\rangle T_{(qN)^\ell\nu}h_{\ell, 0}(\xi)\qquad
\end{align*}
\begin{align*}
&\qquad +\sum_{i=1}^{qN-1}\sum_{\nu \in \Lambda}\Big\langle w_i\big(\sigma-{\mathfrak p}^{-1}N \lambda\big), T_{qN\nu}w_i \Big\rangle T_{(qN)^\ell\nu}h_{\ell, i}(\xi),
\end{align*}

\parindent=0mm \vspace{.0in}
which verifies that $V_{\ell-1}\subset V_\ell \oplus W_\ell $. This completes the proof of Theorem 4.3.\quad\fbox

\parindent=0mm \vspace{.1in}

{\bf{Theorem 4.4.}} {\it For each $\ell,1 \le \ell \le J$, and $i, 0 \le i \le qN-1$, let $w_{\ell,i} \in l^2(\Lambda)$ such that
\begin{align*}
\hat w_{\ell,i}(\xi)= w_{\ell,i0}(\xi)+ \overline{\chi\left(\dfrac{r}{N}\xi\right)}\,w_{\ell,i1}(\xi),\tag{4.4}
\end{align*}

\parindent=0mm \vspace{.0in}
for some ${\mathfrak p}$-periodic functions $ w_{\ell,i0}$ and $w_{\ell,i1}$. For each $\ell$, assume that the matrix ${\mathcal M}_\ell(\xi)$ is unitary, with its entries defined by
\begin{align*}
{\mathcal M}_{\ell, Jt}(\xi)=\dfrac{1}{q\sqrt{ N}}\left\{
\begin{array}{lcr}
 w_{\ell,t0}\left(\xi+\dfrac{ u(J)}{\mathfrak p^{-1}N}\right);\qquad 0\le J \le qN-1;0\le t \le qN-1,\\
 w_{\ell,t1}\left(\xi+\dfrac{u(J-qN)}{\mathfrak p^{-1} N}\right);\qquad qN\le J \le q^2N-1;0\le t \le qN-1,\\
\overline{\chi\left(\dfrac{r}{N}\, \mathfrak p u(J) \right)} w_{\ell,(t-qN)0}\left(\xi+\dfrac{u(J-qN)}{\mathfrak p^{-1} N}\right);\quad 0\le t \le qN-1;\\
\qquad\qquad\qquad \qquad\qquad\qquad\qquad\qquad\qquad qN\le J \le q^2N-1,\\
\overline{\chi\left(\dfrac{r}{N}\, \mathfrak p u(J) \right)} w_{\ell,(t-qN)1}\left(\xi+\dfrac{u(J-qN)}{\mathfrak p^{-1} N}\right);\quad qN\le t,J \le q^2N-1;
\end{array}\right.
\end{align*}
For given $\ell$ and $i$, define $h_{\ell, i}$ as follows
\begin{align*}
\hat h_{\ell, i}(\xi)=\hat h_{\ell-1, i}(\xi)\hat w_{\ell, i}\big(({\mathfrak p}^{-1}N)^{\ell-1}\xi), \quad \text{for}~~ 2 \le \ell \le J
\end{align*}
with $h_{1, i}= w_{1, i}$  and $\hat h_{0, 0}=1$ a.e. Then,
\begin{align*}
V_0=l^2(\Lambda)=V_J\oplus\left(\bigoplus_{m=1}^J W_m\right),\tag{4.5}
\end{align*}
where $W_j=V_{j+1}\ominus V_j, j \in \mathbb Z$ and the $J^{\text{th}}$-stage nonuniform system ${\mathscr H}$ given by (4.1) is orthonormal basis for $l^2(\Lambda)$. }

\parindent=0mm \vspace{.1in}
{\it Proof.} From Theorem 4.2, it follows that  for each $\ell, 1 \le \ell \le J,$ the system
\begin{align*}
\big\{T_{(qN)^\ell\lambda}f_{\ell,i}: 1 \le i \le qN-1, \lambda \in\Lambda\big\}\tag{4.6}
\end{align*}

\parindent=0mm \vspace{.0in}
is  orthonormal in $l^2(\Lambda)$. Therefore, the system $\big\{T_{(qN)^J\lambda}h_{J,0}:  \lambda \in\Lambda\big\}$ and the system defined by (4.1) are both orthonormal for each $\ell$.  Further, using Theorem 4.3, it follows that for each $\ell, 1 \le \ell \le J, V_\ell \subset V_{\ell-1}$ and $V_\ell$ is orthogonal to $W_\ell$ in $V_{\ell-1}$. This means that $W_\ell$ is orthogonal to $W_{\ell-1}$ for each $\ell$. Therefore, the system ${\mathscr H}$ defined by (4.1)is orthonormal in $l^2(\Lambda)$. Since $V_\ell \oplus W_\ell=V_{\ell-1}$, so we can write
\begin{align*}
V_0=V_1\oplus W_1=V_2\oplus W_1\oplus W_2=\cdots=V_J\oplus\left(\bigoplus_{m=1}^J W_m\right).\end{align*}
Since (4.5) holds, the system (4.1) is orthonormal in $l^2(\Lambda)$. This completes the proof of the Theorem 4.4.\quad\fbox

\parindent=0mm \vspace{.1in}

{\bf{Theorem 4.5.}} {\it Under the assumptions of  Theorem 4.4 and for each $\ell \in \mathbb N_0$,  define
\begin{align*}
V_0=l^2(\Lambda),\quad V_\ell=\overline{\text{span}}\big\{T_{(qN)^\ell\lambda}h_{\ell,0}: \lambda \in \Lambda\big\}.
\end{align*}

\parindent=0mm \vspace{.0in}
Then, $\bigcup_{\ell=0}^\infty V_\ell=l^2(\Lambda)$. Also, if $\bigcap_{\ell=0}^\infty V_\ell=\{0\}$, then $l^2(\Lambda)=\bigoplus_{m=1}^\infty W_m$, where $W_j=V_{j+1} \oplus V_j, j \in \mathbb Z$, and  for $\ell \in \mathbb N$, $J^{\text{th}}$-stage nonuniform discrete wavelet system (4.1) is an orthonormal basis for $l^2(\Lambda)$.}

\parindent=0mm \vspace{.2in}
{\it Proof.} Since, for each $\ell \in \mathbb N$, $V_\ell \subset V_{\ell-1}$  and $ V_0=l^2(\Lambda)$, it follows that $\bigcup_{\ell=0}^\infty V_\ell=l^2(\Lambda).$ Using the fact $V_\ell \oplus W_\ell=V_{\ell-1}$, we have
\begin{align*}
V_0=V_1\oplus W_1=V_2\oplus W_1\oplus W_2=\cdots=V_J\oplus\left(\bigoplus_{m=1}^J W_n\right).
\end{align*}

\parindent=0mm \vspace{.0in}
To show $V_0=\bigoplus_{m=1}^\infty W_m$, it is sufficient to show that the orthogonal complement of $\oplus_{m=1}^\infty W_m$ in $V_0$ is $\{0\}$. For this, suppose $f \in V_0$ is orthogonal to $\bigoplus_{m=1}^\infty W_m$. Then $f$ is orthogonal to each $W_m$ for $m \in \mathbb N$. This means that $f$ is a member of each $V_m$ as $W_m$ is orthogonal to $V_m$. Therefore, $f \in \bigcap_{\ell=0}^\infty V_\ell=\{0\}$, which means that $f=0,$ a.e. This completes the proof of Theorem 4.5.\qquad\fbox

\parindent=0mm \vspace{.2in}

{\bf{5.  Connection Between Nonuniform Discrete and Continuous Wavelets}}

\parindent=0mm \vspace{.1in}
In this section, we provide a connection between first-stage nonuniform discrete wavelet system of $l^2(\Lambda )$ and their counterpart nonuniform wavelets  of  $L^2(K)$.

\parindent=0mm \vspace{.1in}
{\bf{Theorem 5.1.}} {\it Let $\big\{\psi_\ell\big\}_{\ell=1}^{qN-1}$ be a system of NUMRA wavelets with scaling function $\psi_0$ in $L^2(K)$. Then, there is a first-stage nonuniform discrete wavelet system for $l^2(\Lambda )$ associated with a system of NUMRA wavelets of  $L^2(K)$.}

\parindent=0mm \vspace{.1in}
{\it{Proof.}} Given a system $\big\{\psi_\ell\big\}_{\ell=1}^{qN-1}$ of NUMRA wavelets with scaling function $\psi_0$ in $L^2(K)$, we define $V_j^*, j \in \mathbb Z$ as $V_j^*=\overline{\text{span}}\big\{D_jT_\lambda \psi_0(x):\lambda \in\Lambda\big\},$  where $\big\{T_\lambda \psi_0(x):\lambda \in\Lambda\big\}$ is an orthonormal basis for $V_0^*$ and the unitary operators $T_\lambda$ and $D_j$ are defined by
\begin{align*}
T_\lambda f(x)= f(x-\lambda), \quad D_j f(x)=(qN)^{j/2}f\big(({\mathfrak p}^{-1}N)^jx\big), \quad \text{for}\, f \in L^2(K).
\end{align*}

\parindent=0mm \vspace{.0in}
Since $\big\{\psi_\ell\big\}_{\ell=1}^{qN-1}\subset V_1^*$, there is  $\big\{w_\ell\big\}_{\ell=1}^{qN-1}\subset l^2(\Lambda )$ such that for each $\ell, 0 \le \ell \le qN-1$,
\begin{align*}
\psi_\ell(x)=\sum_{\lambda \in\Lambda}w_\ell(\lambda)D T_\lambda\psi_0(x). \tag{5.1}
\end{align*}

\parindent=0mm \vspace{.0in}
Equation (5.1) can be written in the frequency domain as
\begin{align*}
\hat \psi_\ell(x)&=\sum_{\lambda \in\Lambda}w_\ell(\lambda)(D T_\lambda\psi_0(x))^{\wedge}=m_\ell\left(\dfrac{\xi}{\mathfrak p^{-1}N}\right)\hat \psi_0\left(\dfrac{\xi}{\mathfrak p^{-1}N}\right),
\end{align*}

\parindent=0mm \vspace{.0in}
where $m_\ell(\xi)=\displaystyle\dfrac{1}{\sqrt{qN}}\sum_{\lambda \in\Lambda}w_\ell(\lambda)\overline{\chi_\lambda(\xi)}$ is $L^2$ locally. Since $\Lambda = \left\{ 0, u(r)/N\right\}+{\mathcal Z}$, we can write
\begin{align*}
m_{\ell}(\xi)= m_{\ell0}(\xi)+ \overline{\chi\left(\dfrac{r}{N}\xi\right)}\,m_{\ell1}(\xi),\quad  0\le \ell\le qN-1, \tag{5.2}
\end{align*}
where $m_{\ell0} $ and $m_{\ell1}$ are locally ${\mathfrak p}$-periodic functions. Therefore, we have  equivalent conditions of orthonormality for the system $\big\{T_\lambda \psi_\ell(x):0\le \ell\le qN-1,\lambda \in\Lambda\big\}$ as
\begin{align*}
\text{(a)}\quad & \sum_{t=0}^{qN-1}\left\{m_{\ell 0} \left(\xi+\dfrac{u(t)}{\mathfrak p^{-1}N}\right)\overline{m_{k 0}\left(\xi+\dfrac{u(t)}{\mathfrak p^{-1}N}\right)}+m_{\ell 0}\left(\xi+\dfrac{u(t)}{\mathfrak p^{-1}N}\right)\overline{m_{k 0}\left(\xi+\dfrac{u(t)}{\mathfrak p^{-1}N}\right)}\right\}= \delta_{\ell, k},\\
\text{(b)}\quad  &\sum_{t=0}^{qN-1}\overline{\chi\left(\dfrac{r}{N}\, \mathfrak p u(t) \right)}\left\{m_{\ell 0}\left(\xi+\dfrac{u(t)}{\mathfrak p^{-1}N}\right)\overline{m_{k 0}\left(\xi+\dfrac{u(t)}{\mathfrak p^{-1}N}\right)}\right.\\
&\qquad\qquad\qquad\qquad\qquad+\left.m_{\ell 0}\left(\xi+\dfrac{u(t)}{\mathfrak p^{-1}N}\right)\overline{m_{k 0}\left(\xi+\dfrac{u(t)}{\mathfrak p^{-1}N}\right)}\right\}=0,
\end{align*}
From the definition of $m_\ell$, we see that
\begin{align*}
m_{\ell}(\xi)= \dfrac{1}{\sqrt{qN}}\hat w_\ell(\xi)
\end{align*}
where $\hat w_\ell$ denotes the Fourier transform in the sense of $l^2(\Lambda )$. Using (3.21), we have
\begin{align*}
m_{\ell0}(\xi)= \dfrac{1}{\sqrt{qN}} w_{\ell0}(\xi)\quad \text{and} \quad m_{\ell1}(\xi)= \dfrac{1}{\sqrt{qN}} w_{\ell1}(\xi),\end{align*}
where $w_{\ell0}$ and $w_{\ell1}$ have same properties as that of $m_{\ell0}$ and $m_{\ell1}$. Substituting the values of $m_{\ell0}$ and $m_{\ell1}$ in (a) and (b), we have for $0\le \ell, k \le qN-1$,
\begin{align*}
\text{(a)}~&\sum_{t=0}^{qN-1}\left\{w_{\ell 0} \left(\xi+\dfrac{u(t)}{\mathfrak p^{-1}N}\right)\overline{w_{k 0}\left(\xi+\dfrac{u(t)}{\mathfrak p^{-1}N}\right)}+w_{\ell 0}\left(\xi+\dfrac{u(t)}{\mathfrak p^{-1}N}\right)\overline{w_{k 0}\left(\xi+\dfrac{u(t)}{\mathfrak p^{-1}N}\right)}\right\}=qN \delta_{\ell, k},\\
\text{(b)}~  &\sum_{t=0}^{qN-1}\overline{\chi\left(\dfrac{r}{N}\, \mathfrak p u(t) \right)}\left\{w_{\ell 0}\left(\xi+\dfrac{u(t)}{\mathfrak p^{-1}N}\right)\overline{w_{k 0}\left(\xi+\dfrac{u(t)}{\mathfrak p^{-1}N}\right)}\right.\\
&\qquad\qquad\qquad\qquad\qquad\qquad\qquad\qquad+\left.w_{\ell 0}\left(\xi+\dfrac{u(t)}{\mathfrak p^{-1}N}\right)\overline{w_{k 0}\left(\xi+\dfrac{u(t)}{\mathfrak p^{-1}N}\right)}\right\}=0.
\end{align*}

\parindent=0mm \vspace{.0in}
These conditions are equivalent to the a.e unitary of the matrix ${\mathcal M}(\xi)$ of order $q^2N \times q^2N$ with entries as in the Corollary 3.8. Therefore,   system ${\mathscr F}(W)$ given by (3.3) is an orthonormal basis of  $l^2(\Lambda )$  and hence, is the first-stage nonuniform discrete wavelet system for $l^2(\Lambda )$. \qquad \fbox

\parindent=8mm \vspace{.1in}

By observing that $m_\ell$ and $w_\ell$ are closely related for each $\ell$, following result can be easily proved

\parindent=0mm \vspace{.1in}
{\bf{Theorem 5.2.}} {\it If ${\mathscr F}(W)$ is the first-stage nonuniform discrete wavelet system as defined by (3.3). Then, there exists a system of NUMRA wavelets $\big\{\psi_\ell\big\}_{\ell=1}^{qN-1}$  with scaling function $\psi_0$ in $L^2(K)$ associated with the first-stage nonuniform discrete wavelet system for $l^2(\Lambda )$.}

\parindent=8mm \vspace{.1in}
It is evident from Theorems 5.1 and 5.2 that the NUMRA wavelets of $L^2(K)$ are connected with the first-stage nonuniform discrete wavelet system of $l^2(\Lambda)$ and vice-versa.

\newpage
\parindent=0mm \vspace{.2in}

{\bf{References}}

\end {document}